\input amstex
\magnification\magstephalf
\documentstyle{amsppt}

\hsize 5.72 truein
\vsize 7.9 truein
\hoffset .39 truein
\voffset .26 truein
\mathsurround 1.67pt
\parindent 20pt
\normalbaselineskip 13.8truept
\normalbaselines
\binoppenalty 10000
\relpenalty 10000
\csname nologo\endcsname 


\font\bc=cmb10
\font\tenbsy=cmbsy10

\catcode`\@=11

\def\myitem#1.{\item"(#1)."\advance\leftskip10pt\ignorespaces}

\def\qedsymbol{{\mathsurround\z@$\square$}}
\redefine\qed{\relaxnext@\ifmmode\let\next\@qed\else
  {\unskip\nobreak\hfil\penalty50\hskip2em\null\nobreak\hfil
    \qedsymbol\parfillskip\z@\finalhyphendemerits0\par}\fi\next}
\def\@qed#1$${\belowdisplayskip\z@\belowdisplayshortskip\z@
  \postdisplaypenalty\@M\relax#1
  $$\par{\lineskip\z@\baselineskip\z@\vbox to\z@{\vss\noindent\qed}}}
\outer\redefine\beginsection#1#2\par{\par\penalty-250\bigskip\vskip\parskip
  \leftline{\tenbsy x\bf#1. #2}\nobreak\smallskip\noindent}
\outer\redefine\genbeginsect#1\par{\par\penalty-250\bigskip\vskip\parskip
  \leftline{\bf#1}\nobreak\smallskip\noindent}

\def\next{\let\@sptoken= }\def\next@{ }\expandafter\next\next@
\def\@futureletnext#1{\let\nextii@#1\futurelet\next\@flti}
\def\@flti{\ifx\next\@sptoken\let\next@\@fltii\else\let\next@\nextii@\fi\next@}
\expandafter\def\expandafter\@fltii\next@{\futurelet\next\@flti}

\let\zeroindent\z@
\let\savedef@\endproclaim\let\endproclaim\relax 
\define\chkproclaim@{\add@missing\endroster\add@missing\enddefinition
  \add@missing\endproclaim
  \envir@stack\endproclaim
  \edef\endit@{\leftskip\the\leftskip\rightskip\the\rightskip}}
\let\endproclaim\savedef@
\def\thing@{.\enspace\egroup\ignorespaces}
\def\thingi@(#1){ \rm(#1)\thing@}
\def\thingii@\cite#1{ \rm\@pcite{#1}\thing@}
\def\thingiii@{\ifx\next(\let\next\thingi@
  \else\ifx\next\cite\let\next\thingii@\else\let\next\thing@\fi\fi\next}
\def\thing#1#2#3{\chkproclaim@
  \ifvmode \medbreak \else \par\nobreak\smallskip \fi
  \noindent\advance\leftskip#1
  \hskip-#1#3\bgroup\bc#2\unskip\@futureletnext\thingiii@}
\let\savedef@\endproclaim\let\endproclaim\relax 
\def\endit{\endproclaim\endit@\let\endit@\undefined}
\let\endproclaim\savedef@

\def\lemma#1{\thing\parindent{Lemma #1}\sl}
\def\prop#1{\thing\parindent{Proposition #1}\sl}
\def\thm#1{\thing\parindent{Theorem #1}\sl}

\def\conj#1{\thing\parindent{Conjecture #1}\sl}

\def\narrowthing#1{\chkproclaim@\medbreak\narrower\noindent
  \it\def\next{#1}\def\next@{}\ifx\next\next@\ignorespaces
  \else\bgroup\bc#1\unskip\let\next\narrowthing@\fi\next}
\def\narrowthing@{\@futureletnext\thingiii@}

\def\@cite#1,#2\end@{{\rm([\bf#1\rm],#2)}}
\def\cite#1{\in@,{#1}\ifin@\def\next{\@cite#1\end@}\else
  \relaxnext@{\rm[\bf#1\rm]}\fi\next}
\def\@pcite#1{\in@,{#1}\ifin@\def\next{\@cite#1\end@}\else
  \relaxnext@{\rm([\bf#1\rm])}\fi\next}

\advance\minaw@ 1.2\ex@
\atdef@[#1]{\ampersand@\let\@hook0\let\@twohead0\brack@i#1,\z@,}
\def\brack@{\z@}
\let\@@hook\brack@
\let\@@twohead\brack@
\def\brack@i#1,{\def\next{#1}\ifx\next\brack@
  \let\next\brack@ii
  \else \expandafter\ifx\csname @@#1\endcsname\brack@
    \expandafter\let\csname @#1\endcsname1\let\next\brack@i
    \else \Err@{Unrecognized option in @[}%
  \fi\fi\next}
\def\brack@ii{\futurelet\next\brack@iii}
\def\brack@iii{\ifx\next>\let\next\brack@gtr
  \else\ifx\next<\let\next\brack@less
    \else\relaxnext@\Err@{Only < or > may be used here}
  \fi\fi\next}
\def\brack@gtr>#1>#2>{\setboxz@h{$\m@th\ssize\;{#1}\;\;$}%
 \setbox@ne\hbox{$\m@th\ssize\;{#2}\;\;$}\setbox\tw@\hbox{$\m@th#2$}%
 \ifCD@\global\bigaw@\minCDaw@\else\global\bigaw@\minaw@\fi
 \ifdim\wdz@>\bigaw@\global\bigaw@\wdz@\fi
 \ifdim\wd@ne>\bigaw@\global\bigaw@\wd@ne\fi
 \ifCD@\enskip\fi
 \mathrel{\mathop{\hbox to\bigaw@{$\ifx\@hook1\lhook\mathrel{\mkern-9mu}\fi
  \setboxz@h{$\displaystyle-\m@th$}\ht\z@\z@
  \displaystyle\m@th\copy\z@\mkern-6mu\cleaders
  \hbox{$\displaystyle\mkern-2mu\box\z@\mkern-2mu$}\hfill
  \mkern-6mu\mathord\ifx\@twohead1\twoheadrightarrow\else\rightarrow\fi$}}%
 \ifdim\wd\tw@>\z@\limits^{#1}_{#2}\else\limits^{#1}\fi}%
 \ifCD@\enskip\fi\ampersand@}
\def\brack@less<#1<#2<{\setboxz@h{$\m@th\ssize\;\;{#1}\;$}%
 \setbox@ne\hbox{$\m@th\ssize\;\;{#2}\;$}\setbox\tw@\hbox{$\m@th#2$}%
 \ifCD@\global\bigaw@\minCDaw@\else\global\bigaw@\minaw@\fi
 \ifdim\wdz@>\bigaw@\global\bigaw@\wdz@\fi
 \ifdim\wd@ne>\bigaw@\global\bigaw@\wd@ne\fi
 \ifCD@\enskip\fi
 \mathrel{\mathop{\hbox to\bigaw@{$%
  \setboxz@h{$\displaystyle-\m@th$}\ht\z@\z@
  \displaystyle\m@th\mathord\ifx\@twohead1\twoheadleftarrow\else\leftarrow\fi
  \mkern-6mu\cleaders
  \hbox{$\displaystyle\mkern-2mu\copy\z@\mkern-2mu$}\hfill
  \mkern-6mu\box\z@\ifx\@hook1\mkern-9mu\rhook\fi$}}%
 \ifdim\wd\tw@>\z@\limits^{#1}_{#2}\else\limits^{#1}\fi}%
 \ifCD@\enskip\fi\ampersand@}

\define\eg{{\it e.g\.}}
\define\ie{{\it i.e\.}}
\define\today{\number\day\ \ifcase\month\or
  January\or February\or March\or April\or May\or June\or
  July\or August\or September\or October\or November\or December\fi
  \ \number\year}
\def\pr@m@s{\ifx'\next\let\nxt\pr@@@s \else\ifx^\next\let\nxt\pr@@@t
  \else\let\nxt\egroup\fi\fi \nxt}

\define\widebar#1{\mathchoice
  {\setbox0\hbox{\mathsurround\z@$\displaystyle{#1}$}\dimen@.1\wd\z@
    \ifdim\wd\z@<.4em\relax \dimen@ -.16em\advance\dimen@.5\wd\z@ \fi
    \ifdim\wd\z@>2.5em\relax \dimen@.25em\relax \fi
    \kern\dimen@ \overline{\kern-\dimen@ \box0\kern-\dimen@}\kern\dimen@}%
  {\setbox0\hbox{\mathsurround\z@$\textstyle{#1}$}\dimen@.1\wd\z@
    \ifdim\wd\z@<.4em\relax \dimen@ -.16em\advance\dimen@.5\wd\z@ \fi
    \ifdim\wd\z@>2.5em\relax \dimen@.25em\relax \fi
    \kern\dimen@ \overline{\kern-\dimen@ \box0\kern-\dimen@}\kern\dimen@}%
  {\setbox0\hbox{\mathsurround\z@$\scriptstyle{#1}$}\dimen@.1\wd\z@
    \ifdim\wd\z@<.28em\relax \dimen@ -.112em\advance\dimen@.5\wd\z@ \fi
    \ifdim\wd\z@>1.75em\relax \dimen@.175em\relax \fi
    \kern\dimen@ \overline{\kern-\dimen@ \box0\kern-\dimen@}\kern\dimen@}%
  {\setbox0\hbox{\mathsurround\z@$\scriptscriptstyle{#1}$}\dimen@.1\wd\z@
    \ifdim\wd\z@<.2em\relax \dimen@ -.08em\advance\dimen@.5\wd\z@ \fi
    \ifdim\wd\z@>1.25em\relax \dimen@.125em\relax \fi
    \kern\dimen@ \overline{\kern-\dimen@ \box0\kern-\dimen@}\kern\dimen@}%
  }

\catcode`\@\active

\let\PVstyle=d 
 
\font\tenscr=rsfs10 
\font\sevenscr=rsfs7 
\font\fivescr=rsfs5 
\skewchar\tenscr='177 \skewchar\sevenscr='177 \skewchar\fivescr='177
\newfam\scrfam \textfont\scrfam=\tenscr \scriptfont\scrfam=\sevenscr
\scriptscriptfont\scrfam=\fivescr
\define\scr#1{{\fam\scrfam#1}}
\let\Cal\scr
 
\let\0\relax 
\mathchardef\idot="202E
\define\restrictedto#1{\big|_{#1}}
\define\na{{\text{na}}}
\define\ord{\operatorname{ord}}
\define\red{{\text{red}}}
\define\Spec{\operatorname{Spec}}
\define\Supp{\operatorname{Supp}}
 
\topmatter
\title A more general abc conjecture\endtitle
\author Paul Vojta\endauthor
\affil University of California, Berkeley\endaffil
\address University of California, Department of Mathematics \#3840,
  Berkeley, CA \ 94720--3840, USA\endaddress
\date 1 May 1998\enddate
\thanks Supported by NSF grant DMS95-32018, the Institute for Advanced Study,
  and IHES\endthanks

\abstract
This note formulates a conjecture generalizing both the $abc$ conjecture of
Masser-Oesterl\'e and the author's diophantine conjecture for algebraic
points of bounded degree. It also shows that the new conjecture is implied
by the earlier (apparently weaker) conjecture. 

As with most of the author's conjectures, this new conjecture stems from
analogies with Nevanlinna theory; in this case it corresponds to a
Second Main Theorem in Nevanlinna theory with truncated counting functions.
The original abc conjecture of Masser and Oesterl\'e corresponds to the
Second Main Theorem with truncated counting functions on $\Bbb P^1$
for the divisor $[0]+[1]+[\infty]$. 
\endabstract
\endtopmatter
 
\document

This paper appeared in the {\it International Mathematics Research Notices,}
{\bf 1998} (1998) 1103--1116.

\bigbreak
 
In this note we formulate a conjecture generalizing both the $abc$ conjecture
of Masser-Oesterl\'e and the author's diophantine conjecture for algebraic
points of bounded degree.  We also show that the earlier (apparently weaker)
conjecture implies the new conjecture.

As with most of the author's conjectures, this new conjecture stems from
analogies with Nevanlinna theory.  In this particular case the conjecture
corresponds to relacing the usual counting function of Nevanlinna theory
with a truncated counting function.  In particular, the $abc$ conjecture
of Masser and Oesterl\'e corresponds to Nevanlinna's Second Main Theorem
with truncated counting functions applied to the divisor $[0]+[1]+[\infty]$
on $\Bbb P^1$.

The first section of this paper introduces the notation that will be
used throughout the paper.  Section \02 formulates the new conjecture
and discusses some examples related to the new conjecture, including
an ``$abcde\dots$ conjecture'' and a conjecture of Buium.
The third and final section of this paper shows that the new conjecture
is implied by the (apparently weaker) older conjecture without truncated
counting functions.

\beginsection{\01}{Notation}

This section briefly recalls the Nevanlinna-based notation from
\cite{V, \S3.2} that will be needed for stating the conjecture.

Let $k$ be a global field.  Its set of places will be denoted $M_k$.
Each place $v\in M_k$ has an associated almost absolute value $\|\cdot\|_v$,
normalized as follows.  If $k$ is a number field, then let $\Cal O_k$
denote its ring of integers.  A non-archimedean place $v\in M_k$ corresponds
to a nonzero prime ideal $\frak p\subseteq\Cal O_k$, and we set
$\|x\|_v=(\Cal O_k:\frak p)^{-\ord_{\frak p}x}$ if $x\in k^\times$
(and $\|0\|_v=0$).  If $v$ is archimedean, then $v$ corresponds to a real
embedding $\sigma\:k\hookrightarrow\Bbb R$ or a complex conjugate pair
of complex embeddings $\sigma,\bar\sigma\:k\hookrightarrow\Bbb C$,
and we set $\|x\|_v=|\sigma(x)|$ or $\|x\|_v=|\sigma(x)|^2$, respectively.
(In the latter case the triangle inequality fails to hold, hence the
terminology ``almost absolute value.'')  If $k$ is a function field,
then $k=K(C)$ for some smooth projective curve $C=C_k$ over the field
of constants $k_0\subseteq k$; places $v\in M_k$ correspond bijectively to
closed points $p\in C_k$, and we set $\|x\|_v=\exp(-[K(p):k_0]\ord_p(x))$
if $x\in k^\times$.

If $k$ is a global field and $v\in M_k$, then let $k_v$ denote the completion
of $k$ at $v$, and let $\Bbb C_v$ denote the completion of the algebraic closure
of $k_v$.

The proofs in this paper will use Arakelov theory, some of whose language
is as follows.  The general idea is that the number field case should mimic
as closely as possible the situation encountered in the function field case,
where one has the projective curve $C_k$, and one can work with intersection
theory on a proper scheme over $C_k$.  In the number field case this is
accomplished by formally adding analytic information for each archimedean
place; hence the role of $C_k$ is played by an arithmetic scheme $\Bbb M_k$
consisting of $\Spec\Cal O_k$, with finitely many points added,
corresponding to the archimedean places.  Therefore, one can think of
$\Bbb M_k$ as an object whose closed points are in canonical bijection
with $M_k$.  We also define the {\bc non-archimedean part} $(\Bbb M_k)_\na$
of $\Bbb M_k$ to be $\Spec\Cal O_k$.  In the function field case,
we set $\Bbb M_k=(\Bbb M_k)_\na=C_k$.

In this paper, a {\bc variety} over a field $k$ is an integral separated
scheme of finite type over $k$.  If $k$ is a global field (which we assume from
now on), then an {\bc arithmetic variety} $\Cal X$ over $\Bbb M_k$ is an
integral scheme $\Cal X_\na$, flat, separated, and of finite type over
$(\Bbb M_k)_\na$, together with some analytic information at the archimedean
places (which will play no role in this paper).  An arithmetic variety
$\Cal X$ over $\Bbb M_k$ is {\bc proper} over $\Bbb M_k$ if $\Cal X_\na$
is proper over $(\Bbb M_k)_\na$.

If $X$ is a variety over $k$, then a {\bc model} for $X$ is an arithmetic
variety $\Cal X$ over $\Bbb M_k$, together with an isomorphism 
$X\cong\Cal X\times_{\Bbb M_k} k:=\Cal X_\na\times_{(\Bbb M_k)_\na} k$.

Let $X$ be a complete variety over a global field $k$, and let $\Cal X$
be a proper model over $\Bbb M_k$.  An algebraic point $P\in X(\bar k)$
determines a map $\sigma\:\Bbb M_E\to\Cal X$ over $\Bbb M_k$
(that is, a map $\sigma_\na\:(\Bbb M_E)_\na\to\Cal X_\na$
over $(\Bbb M_k)_\na$), where $E=k(P)$.

A {\bc Cartier divisor} $D$ on an arithmetic variety $\Cal X$ over $\Bbb M_k$
is a Cartier divisor $D_\na$ on $\Cal X_\na$, together with
Green functions $g_{D,v}$ for $D$ on $\Cal X(\Bbb C_v):=\Cal X_\na(\Bbb C_v)$
for all archimedean $v\in M_k$.  The Green functions are taken to be normalized
corresponding to $-\log\|\cdot\|_v$.  A {\bc principal} Cartier divisor may
be defined in the obvious way, using $g_{(f),v}=-\log\|f\|_v$, and one
obtains the notion of linear equivalence of Cartier divisors on $\Cal X$.
A {\bc line sheaf} $\Cal L$ on $\Cal X$ is a line sheaf $\Cal L_\na$ on
$\Cal X_\na$, together with metrics at the archimedean places.  We have
a natural bijection between the group of divisor classes and the group of
isomorphism classes of line sheaves on $\Cal X$.

One may regard $\Bbb M_k$ as an arithmetic variety over itself.  Let $D$
be a Cartier divisor on $\Bbb M_k$, and let $v\in M_k$.  Then we define
the {\bc degree} of $D$ at $v$ as follows.  If $v$ is archimedean,
then $\Bbb M_k(\Bbb C_v)$ consists of just one point, so $g_{D,v}$ is just
a real number, and we let $\deg_v D=g_{D,v}$.  Otherwise, $v$ corresponds
to a closed point on $(\Bbb M_k)_\na$, also denoted $v$.  Let $n_v$ be the
multiplicity of $D_\na$ at that point, let $K(v)$ denote the residue field
at $v$, and let
$$\mu(K(v))
  = \cases \log\#K(v)&\qquad\text{if $k$ is a number field; or} \\
    [K(v):k_0]&\qquad\text{if $k$ is a function field with
      field of constants $k_0$.}\endcases$$
We then define
$$\deg_v D = n_v\mu(K(v))\;.$$
For subsets $S\subseteq M_k$, let
$$\deg_S D = \sum_{v\in S}\deg_v D\;.$$
For $S=M_k$, we define $\deg D=\deg_{M_k} D$.  If $f\in k^\times$, then
the Artin-Whaples product formula implies that $\deg(f)=0$.
Thus the degree $\deg\Cal L$ of a line sheaf $\Cal L$ on $\Bbb M_k$ is
well defined, via the corresponding divisor class.

If $E$ is a finite extension of $k$, then $\Bbb M_E$ may be regarded as
an arithmetic variety over $\Bbb M_k$, as well as an arithmetic variety
over itself.  The definitions of Cartier divisors and line sheaves
on $\Bbb M_E$ do not coincide in this case (due to possible archimedean
places), but there are obvious translations back and forth.

Let $k$ be a global field, let $X$ be a complete variety over $k$,
let $\Cal X$ be a proper model for $X$ over $\Bbb M_k$, let $\Cal L$ be a
line sheaf on $\Cal X$, let $P\in X(\bar k)$, let $E$ be a finite extension
of $k$ containing $k(P)$, and let $\sigma\:\Bbb M_E\to\Cal X$
correspond to $P$.  Then the height of $P$ (relative to $\Cal L$ and $k$)
satisfies
$$h_{\Cal L,k}(P) = \frac1{[E:k]}\deg\sigma^{*}\Cal L\;.$$
It is linear and functorial in $\Cal L$, and is independent of
the choice of $E$.  If $D$ is a Cartier divisor on $\Cal X$ such that
$P\notin\Supp D$, then $\sigma^{*}D$ is defined, and we have
$$h_{\Cal O(D),k}(P) = \frac1{[E:k]}\deg\sigma^{*}D\;.$$
Let $S$ be a finite set of places of $\Bbb M_k$ containing the archimedean
places, and let $T=\{w\in M_E\bigm|\text{$w\mid v$ for some $v\in S$}\}$.
Then we may split $h_{\Cal O(D),k}(P)$ into two terms:
$$h_{\Cal O(D),k}(P) = m_{k,S}(D,P) + N_{k,S}(D,P)\;,$$
where
$$m_{k,S}(D,P) := \frac1{[E:k]}\deg_T\sigma^{*}D
  \qquad\text{and}\qquad
  N_{k,S}(D,P) := \frac1{[E:k]}\deg_{M_E\setminus T}\sigma^{*}D$$
are called the {\bc proximity function} and {\bc counting function},
respectively.  (The names and notation come from Nevanlinna theory.)

We also define the notation $w\mid S$ to mean $w\in T$, where $T$ is defined
above.

If $\Cal X'$ is another model for $\Cal X$, and if $D'$ is a Cartier divisor
on $\Cal X'$ coinciding on $X$ with $D$, then we have
$$\align h_{\Cal O(D'),k}(P) &= h_{\Cal O(D),k}(P) + O(1)\;,\\
  m_{k,S}(D',P) &= m_{k,S}(D,P) + O(1)\;,\qquad\text{and} \\
  N_{k,S}(D',P) &= N_{k,S}(D,P) + O(1)\endalign$$
for all $P\in X(\bar k)$.  Therefore the height, proximity, and counting
functions may be discussed in terms of Cartier divisors on $X$ if their
values are only needed up to a bounded function.

For more details on height, proximity, and counting functions, including
an alternative definition using Weil functions, see \cite{V, \S3.4}.

Since $k$ and $S$ will often be fixed for a given discussion, they will
often be omitted from the notation.

All places in $M_k\setminus S$ are non-archimedean; hence the counting function
may be written
$$N_{k,S}(D,P)
  = \frac1{[E:k]}\sum\Sb w\in M_E\\w\nmid S\endSb n_v\cdot\mu(K(v))\;.$$
Here $n_v$ is the multiplicity of $v$ in the divisor $\sigma^{*}D$, as above,
and $\mu(K(v))$ is as defined earlier.

We may then define the {\bc truncated counting function}
$$N_{k,S}^{(1)}(D,P)
  = \frac1{[k(P):k]}\sum\Sb w\in M_{k(P)}\\w\nmid S\endSb
    \min\{1,n_v\}\cdot\mu(K(v))$$
for effective Cartier divisors $D$ on $\Cal X$ (\ie, Cartier divisors such
that $D_\na$ is effective) and for points $P\in X(\bar k)\setminus\Supp D$.
The truncated counting function is not necessarily additive or functorial
in $D$.

We next define some quantities related to the logarithm of the
discriminant of a number field, or the genus of the curve corresponding
to a function field.  Let $E$ be a finite extension of a global field $k$
of characteristic $0$.  Then we have a finite morphism $\Bbb M_E\to\Bbb M_k$.
In this situation, let $R_{E/k}$ denote the Cartier divisor on $\Bbb M_E$
such that $(R_{E/k})_\na$ is the ramification divisor of the corresponding
map $(\Bbb M_E)_\na\to(\Bbb M_k)_\na$, and such that the corresponding
Green functions are all zero.  We then define
$$d_k(E) = \frac1{[E:k]}\deg R_{E/k}
  \qquad\text{and}\qquad
  d_{k,S}(E) = \frac1{[E:k]}\deg_{M_E\setminus S} R_{E/k}\;,$$
where $S$ is a finite subset of $M_k$ containing the archimedean places,
and $M_E\setminus S$ means $\{w\in M_E\bigm|w\nmid S\}$.  We note that:
\roster
\myitem i.  If $k$ is a number field and $D_k$ denotes its discriminant, then
$$d_k(E) = \frac1{[E:k]}\log|D_E| - \log|D_k|\;;$$
\myitem ii.
$$0 \le d_k(E) - d_{k,S}(E) \le O_{[E:k],S}(1)\;;\tag\01.1$$
and
\myitem iii.  if $F$ is a finite extension of $E$, then
$$d_{k,S}(F) - d_{k,S}(E)
  = \frac1{[F:k]}\deg_{M_F\setminus S} R_{F/E}\;.\tag\01.2$$
\endroster
For $P\in X(\bar k)$, we define $d_k(P)=d_k(k(P))$
and $d_{k,S}(P)=d_{k,S}(k(P))$.

Let $X$ be a smooth variety.  A {\bc normal crossings divisor}
on $X$ is a divisor that, for all points $P\in X$, can be represented
in the completed local ring $\widehat{\Cal O}_{P,X}$ by a principal divisor
$(x_1\dotsm x_r)$, where $x_1,\dots,x_r$ form a part of a regular sequence
for $\widehat{\Cal O}_{P,X}$.  (A normal crossings divisor must therefore
be effective, and all irreducible components of its support must occur
with multiplicity $1$.)  We also say that a divisor $D$ {\bc has normal
crossings} if it is a normal crossings divisor.

Finally, a divisor $D$ or line sheaf $\Cal L$ on a complete variety $X$
is said to be {\bc big} if there is a constant $c>0$ such that
$$h^0(X,\Cal O(nD))\ge cn^{\dim X}
  \qquad\text{or}\qquad
  h^0(X,\Cal L^{\otimes n})\ge cn^{\dim X}\;,$$
respectively, for all sufficiently large and divisible integers $n$.
By Kodaira's lemma (see, for example, \cite{V, Prop.~1.2.7}), if $D$ is a big
divisor and $A$ is an ample divisor, then $nD-A$ is linearly equivalent
to an effective divisor for sufficiently large and divisible integers $n$.

\beginsection{\02}{The conjecture}

We begin by recalling from \cite{V, Conj.~5.2.6} the general conjecture on
algebraic points of bounded degree.

\bgroup\narrower\it
Throughout this section, $k$ is a global field of characteristic $0$
and $S$ is a finite set of places of $k$ containing the archimedean places.
\par\egroup

\conj{\02.1}  Let $X$ be a smooth complete variety over $k$, let $D$ be a
normal crossings divisor on $X$, let $\Cal K$ denote the canonical line sheaf
on $X$, let $\Cal A$ be a big line sheaf on $X$, let $r\in\Bbb Z_{>0}$,
and let $\epsilon>0$.  Then there exists a proper Zariski-closed subset
$Z=Z(k,S,X,D,\Cal A,r,\epsilon)\subsetneqq X$ such that
$$h_{\Cal K}(P) + m(D,P) \le d_k(P) + \epsilon\,h_{\Cal A}(P) + O(1)
  \tag\02.1.1$$
for all $P\in X(\bar k)\setminus Z$ with $[k(P):k]\le r$.
\endit

(In \cite{V}, the term $d_k(P)$ had a factor $\dim X$ in front, but in
recent years it has become apparent that the inequality may be true without
this factor.)

Since $h_{\Cal O(D)}(P)=m(D,P)+N(D,P)+O(1)$, (\02.1.1) is equivalent to
$$N(D,P) + d_k(P) \ge h_{\Cal K(D)}(P) - \epsilon\,h_{\Cal A}(P) - O(1)\;.
  \tag\02.2$$
One may then ask whether $N(D,P)$ could be replaced by the truncated counting
function:

\conj{\02.3}  Conjecture \02.1 holds with (\02.1.1) replaced by
$$N^{(1)}(D,P) + d_k(P)
  \ge h_{\Cal K(D)}(P) - \epsilon\,h_{\Cal A}(P) - O(1)\;.\tag\02.3.1$$
\endit

We always have $N^{(1)}(D,P)\le N(D,P)$, so Conjecture \02.3 is obviously
stronger than Conjecture \02.1.  The main goal of the next section will be to
show that the converse holds; \ie, Conjecture \02.1 implies Conjecture \02.3.

First, however, we shall show how Conjecture \02.3 generalizes the
``$abc$ conjecture'' of Masser and Oesterl\'e, which is the following:

\conj{\02.4} (Masser-Oesterl\'e)  For all $\epsilon>0$ there is a constant
$C>0$ such that for all $a,b,c\in\Bbb Z$ with $a+b+c=0$ and $(a,b,c)=1$,
we have
$$\max\{|a|,|b|,|c|\} \le C\cdot\prod_{p\mid abc} p^{1+\epsilon}\;.\tag\02.4.1$$
\endit

To relate this conjecture to Conjecture \02.3, we begin by translating
(\02.4.1) into the language of Section \01.  The triple $(a,b,c)$ determines
a point $P:=[a:b:c]$ on $\Bbb P^2$, which in fact lies on the line
$x_0+x_1+x_2=0$.  The height of this point is
$h_{\Cal O(1),\Bbb Q}(P) = h(P) = \log\max\{|a|,|b|,|c|\}$ (since $(a,b,c)=1$).
The relative primeness condition also implies that the curve on
$\Bbb P^2_{\Bbb Z}$ corresponding to $P$ meets the divisor $[x_0=0]$
at a prime $p$ if and only if $p\mid a$.  Similarly, it meets the divisors
$[x_1=0]$ and $[x_2=0]$ at $p$ if and only if $p\mid b$ and $p\mid c$,
respectively.  Let $D$ be the divisor $[x_0=0]+[x_1=0]+[x_2=0]$ on $\Bbb P^2$.
Then
$$N_{\Bbb Q,\{\infty\}}^{(1)}(D,P) = \sum_{p\mid abc}\log p\;.$$
Thus (\02.4.1) can be written
$$h(P) \le (1+\epsilon)N^{(1)}(D,P) + O_\epsilon(1)$$
or (with a different $\epsilon$)
$$(1-\epsilon)h(P) \le N^{(1)}(D,P) + O_\epsilon(1)\tag\02.5$$
for all $P\in\Bbb P^2(\Bbb Q)$ lying on the line $x_0+x_1+x_2=0$.
Let $X$ be this line.  Its canonical line sheaf is $\Cal K\cong\Cal O(-2)$,
and $D\restrictedto X$ consists of three distinct points, so
$\Cal K(D)\cong\Cal O(1)$; hence $h(P)=h_{\Cal K(D)}(P) + O(1)$.
Also let $\Cal A=\Cal O(1)$; then (\02.5) becomes
$$N^{(1)}(D,P)
  \ge h_{\Cal K(D)}(P) - \epsilon\,h_{\Cal A}(P) - O_\epsilon(1)\;,$$
which coincides with (\02.3.1) since we are dealing with rational points
and therefore $d(P)=0$ for all $P$.

Thus, it follows that the $abc$ conjecture coincides with the special case
of Conjecture \02.3 when $k=\Bbb Q$, $S=\{\infty\}$, $X=\Bbb P^1$, $D$
consists of three distinct points, and $r=1$.  Conjecture \02.3 can
therefore be viewed as doing for the $abc$ conjecture what Conjecture \02.1
did for Roth's theorem.

One may wonder what this says about what the exponent should be for the $abc$
conjecture in more than three variables (e.g., $a+b+c+d=0$).
Conjecture \02.3 suggests that the exponent should still be $1+\epsilon$,
but only generically.  Indeed, given $n\ge3$ let $X$ be the hyperplane
$x_0+\dots+x_{n-1}=0$ in $\Bbb P^{n-1}$, and let $D$ be the restriction to
$X$ of the sum of the coordinate hyperplanes.  Then $X$ is smooth and $D$ has
normal crossings.  If $x_0,\dots,x_{n-1}$ are nonzero integers satisfying
$x_0+\dots+x_{n-1}=0$ and $(x_0,\dots,x_{n-1})=1$, then
$P:=[x_0:\dots:x_{n-1}]\in X$, $h(P)=\log\max\{|x_0|,\dots,|x_{n-1}|\}$,
and $N^{(1)}(D,P)=\log\prod_{p\mid x_0\dotsm x_{n-1}}p$.
Then Conjecture \02.3 would imply that
$$\max\{|x_0|,\dots,|x_{n-1}|\}
  \le C\prod_{p\mid x_0\dotsm x_{n-1}}p^{1+\epsilon}$$
for all $x_0,\dots,x_{n-1}$ as above {\it outside a proper Zariski-closed
subset.}

This subset is, in fact, essential:  consider the following example.
It is well known that there exist infinitely many triples $(a,b,c)\in\Bbb Z^3$
with $a+b+c=0$, $(a,b,c)=1$, and $\max\{|a|,|b|,|c|\}\ge\prod_{p\mid abc}p$.
For such triples, we have
$$a^2+2ab+b^2-c^2=0\;.$$
Letting $a'=a^2$, $b'=2ab$, $c'=b^2$, and $d'=-c^2$, we have
$$\align a'+b'+c'+d' &= 0 \;, \\
  (a',b',c',d') &= 1 \;, \\
  h([a':b':c':d']) &= 2h([a:b:c])+O(1) \;\text{, and} \\
  \prod_{p\mid a'b'c'd'}p &= \prod_{p\mid abc}p\;.\endalign$$
Thus there are infinitely many points with
$$\max\{|a'|,|b'|,|c'|,|d'|\} \gg \prod_{p\mid a'b'c'd'}p^2\;.$$

This does not contradict Conjecture \02.3, however, because of the exceptional
subset $Z$.  Instead, it shows that working with $Z$ is the hardest part
in determining what the conjectural exponent should be.

Even with an exceptional subset, and requiring the $x_i$ to be pairwise
relatively prime, though, it is fairly easy to see that an exponent better
than $1$ is not possible:

\prop{\02.6}  Let $n\ge 3$ be an integer and let $\epsilon>0$ be arbitrary.
Then there exists a Zariski-dense set of points on the hyperplane
$x_0+\dots+x_{n-1}$ in $\Bbb P^{n-1}$ having pairwise relatively prime
homogeneous coordinates $[x_0:\dots:x_{n-1}]$ such that
$$\prod_{p\mid x_0\dotsm x_{n-1}}p
  \le \epsilon\max\{|x_0|,\dots,|x_{n-1}|\}\;.$$
\endit

\demo{Proof}  For $n=3$ this is already well known.

Assume for now that $n=4$.  Since $\log_9 25$ is irrational,
its positive integer multiples are dense in $\Bbb R/2\Bbb Z$; hence
there exist positive integers $e_1$ and $e_2$ with $e_1$ odd, such
that $0<9^{e_1}-25^{e_2}<\epsilon\cdot9^{e_1}$.  Let $x_1=-9^{e_1}$,
$x_2=25^{e_2}$, $x_3=1$, and choose $x_0$ so that the sum vanishes.
This gives a tuple $(x_0,x_1,x_2,x_3)$ whose elements are pairwise relatively
prime, whose elements add up to zero, and which satisfies
$$\prod_{p\mid x_0x_1x_2x_3}p \le 15|x_0|
  \le 15\epsilon\max\{|x_0|,|x_1|,|x_2|,|x_3|\}\;.$$
The set of all such points is Zariski-dense in the hyperplane
$x_0+x_1+x_2+x_3=0$:  otherwise, some irreducible curve would contain
infinitely many of them, and hence some irreducible polynomial $f(X,Y)$
would satisfy $f(9^{e_1},25^{e_2})=0$ for infinitely many pairs $(e_1,e_2)$.
Applying the unit equation for $\Bbb Z[1/15]$ to the terms of such equations
gives finitely many linear relations in those terms, leading to a
contradiction.

Finally, assume that $n\ge 5$.  Let $p_1,\dots,p_{n-1}$ be distinct primes
greater than $n$, let $r_1,\dots,r_{n-1}$ be positive integers such that
$p_i^{r_i}\equiv1\pmod{p_j}$ for all $i\ne j$, and let $q_i=p_i^{r_i}$
for all $i$.  As before, there exist infinitely many pairs $(e_1,e_2)$ of
positive integers such that
$$0 < q_1^{e_1}-q_2^{e_2} < \epsilon q_1^{e_1}\;.$$
For such pairs, let $e_3,\dots,e_{n-1}$ be nonnegative integers such that
$$q_3^{e_3}+\dots+q_{n-1}^{e_{n-1}} < q_1^{e_1}-q_2^{e_2}\;,$$
let $x_0$ be the difference, let $x_1=-q_1^{e_1}$, and let $x_i=q_i^{e_i}$
for $i=2,\dots,n-1$.  Then one will have $x_0+\dots+x_{n-1}=0$,
the $x_i$ will be pairwise relatively prime, and the inequality
$$\prod_{p\mid x_0\dotsm x_{n-1}}p \le p_1\dotsm p_{n-1}|x_0|
  \le p_1\dotsm p_{n-1}\epsilon\max\{|x_0|,\dots,|x_{n-1}|\}$$
will hold.  Moreover, these points will lie outside any
given proper Zariski-closed subset:  if $e_1$ and $e_2$ are large, then
$q_1^{e_1}-q_2^{e_2}$ will be large (\eg, by Baker's theorem), and then
there will be enough choices for each of the remaining $e_i$ to avoid
the subset.\qed
\enddemo

Finally, we mention a conjecture of A. Buium \cite{B}:

\conj{\02.7}  Let $A$ be an abelian variety over $k$, let $D$ be an
ample effective divisor on $A$, and let $\Cal A$ be an ample line sheaf
on $A$.  Then $N^{(1)}(D,P) \gg h_D(P) + O(1)$ for all $P\in A(k)\setminus D$.
\endit

This, too, follows from Conjecture \02.3.  Indeed, let $X$ be a closed
subvariety of $A$ not contained in the support of $D$.  There exists a proper
birational morphism $\pi\:X'\to X$ such that $X'$ is nonsingular and
$D':=(\pi^{*}D)_\red$ is a normal crossings divisor on $X'$.
Rational points $P\in X(k)$ lying outside a proper Zariski-closed subset $Z$
may be lifted to $P'\in X'(k)$; for these points Conjecture \02.3 gives
$$\split N^{(1)}(D,P) &= N^{(1)}(D',P') + O(1) \\
  &\ge h_{\Cal K'(D')}(P') - \epsilon\,h_{\Cal A'}(P') - O(1)\endsplit$$
(after enlarging $Z$).  Here $\Cal K'$ is the canonical line sheaf on $X'$
and $\Cal A'$ is a big line sheaf on $X'$.  But $X$ has Kodaira
dimension $\ge0$ and $D$ is ample, so $\Cal K'(D')$ is big; hence (after
enlarging $Z$ again) we have
$$h_{\Cal K'(D')}(P') - \epsilon\,h_{\Cal A'}(P') \gg h_D(P) + O(1)$$
for suitable $\epsilon>0$.  Starting with $X=A$ and repeating with smaller
and smaller $X$ coming from the exceptional Zariski-closed subset, we then
conclude by noetherian induction that Conjecture \02.7 follows from
Conjecture \02.3.

\beginsection{\03}{The harder implication}

The goal of this section is to prove:

\thm{\03.1}  Conjecture \02.1 implies Conjecture \02.3.
\endit

The general strategy of the proof is to start with the data $X$, $D$, etc.
for which one wants to derive (\02.3.1), determine a large integer $e$
depending on these data, and construct a cover $X'$ of $X$, ramified to
order $\ge e$ everywhere above $D$ and unramified elsewhere.  Then a point
$P\in X(\bar k)$ of bounded degree will lift to a point $P'\in X'(\bar k)$
of bounded (but larger) degree, and the ramification of $k(P')$ over $k(P)$
will occur only at places contributing to $N(D,P)$.  But the contribution
to $d(P')-d(P)$ at a place $w$ of $k(P')$ is limited, so in fact $d(P')-d(P)$
is more closely related to $N^{(1)}(D,P)$.  One can then apply Conjecture \02.1
on $X'$ to deduce Conjecture \02.3 for points $P$ on $X$ outside a
proper Zariski-closed subset.

This general method of proof has been used previously by the author
\cite{V, pp.~71--72}, and by Darmon and Granville \cite{D-G}.

We start with a lemma describing how $K_X+D$ changes when pulled back via a
generically finite morphism.

\lemma{\03.2}  Let $\pi\:X'\to X$ be a generically finite morphism of
smooth varieties, and let $D$ and $D'$ be normal crossings divisors on
$X$ and $X'$, respectively, such that $\Supp D'=(\pi^{*}D)_\red$.
Let $\Cal K$ and $\Cal K'$ be the canonical line sheaves on $X$ and $X'$,
respectively.  Then
$$\Cal K'(D') \ge \pi^{*}(\Cal K(D))$$
relative to the cone of line sheaves with $h^0>0$.  Moreover,
$\Cal K'(D')\otimes\pi^{*}(\Cal K(D))\spcheck$ has a global section
vanishing only on the support of $D'$ or where $\pi$ ramifies.
\endit

\demo{Proof}  We have a natural map
$\pi^{*}\Omega^1_X[\log D]\to\Omega^1_{X'}[\log D']$ which is an isomorphism
at generic points of $X'$; hence taking the maximal exterior power gives
an injection of sheaves $\pi^{*}(\Cal K(D))\to\Cal K'(D')$.
See also \cite{I, \S11.4a}.\qed
\enddemo

Consequently, one may use Chow's lemma and resolution of singularities to
find a smooth projective variety $X'$ and a proper birational morphism
$\pi\:X'\to X$ such that $D':=(\pi^{*}D)_\red$ has normal crossings;
the lemma then shows that
$$h_{\Cal K'(D')}\ge h_{\Cal K(D)}\circ\pi+O(1)$$
outside of a proper Zariski-closed subset.  Since $\Supp D'=\pi^{-1}(\Supp D)$,
we have $N^{(1)}(D',P)=N^{(1)}(D,\pi(P))+O(1)$ for $P\in X'(\bar k)$;
hence Conjecture \02.3 for $X'$ and $D'$ implies the same conjecture
for $X$ and $D$.  Thus we may assume that $X$ is projective.

Next, we note that if $D$ and $D_1$ are effective divisors such that $D+D_1$
has normal crossings, then Conjecture \02.3 for $D+D_1$ implies that the
conjecture holds also for $D$.  Indeed,
$$\split N^{(1)}(D+D_1,P) - N^{(1)}(D,P) &\le N^{(1)}(D_1,P) \\
  &\le N(D_1,P) \\
  &\le h_{\Cal O(D_1)}(P) + O(1) \;,\endsplit$$
so (\02.3.1) with $D$ replaced by $D+D_1$ implies the original (\02.3.1).

Thus, we may assume that $D$ is very ample; in addition, we will enlarge $D$
at one point in the proof.

\lemma{\03.3}  Let $e\in\Bbb Z_{>0}$ and let $D_1$ be an effective divisor
such that $D_1\sim D$ and $D+D_1$ has normal crossings.  Then there is a
smooth variety $X_1$ and a proper generically finite morphism $\pi_1\:X_1\to X$
such that the support of the ramification divisor of $\pi_1$ is equal to
$\pi_1^{-1}(\Supp(D+D_1))$, all components of $\pi^{*}(D+D_1)$ have
multiplicity $\ge e$ unless they lie over $D\cap D_1$,
$(\pi^{*}(D+D_1))_\red$ has normal crossings,
and $K(X_1)=K(X)(\root e\of f)$ for some $f\in K(X)^{*}$.
\endit

\demo{Proof}  Pick $f\in K(X)^{*}$ such that $(f)=D-D_1$.
Let $\pi_0\:X_0\to X$ be the normalization of $X$ in the field
$K(X)(\root e\of f)$, let $\pi\:X_1\to X_0$ be a desingularization of $X_0$
such that $\bigl(\pi^{*}(\pi_0^{*}(D+D_1))\bigr)_\red$ has normal crossings,
and let $\pi_1=\pi_0\circ\pi$.  We may assume that $\pi_1$ is \'etale
outside of $\pi_0^{-1}(\Supp D)$.  If $P\notin D\cap D_1$, then there is an
open neighborhood $U$ of $P$ such that $D\restrictedto U=(f)$ and
$D_1\restrictedto U=0$, or $D\restrictedto U=0$ and $D_1\restrictedto U=(1/f)$.
Then $\pi_0^{*}(D+D_1)$ over $U$ is $e$ times the (principal) divisor
$(\root e\of f)$ or $(1/\root e\of f)$.  Therefore all components
of $\pi_1^{*}(D+D_1)$ meeting $\pi_1^{-1}(U)$ have multiplicity divisible
by $e$.\qed
\enddemo

\lemma{\03.4}  There exists a normal crossings divisor $D^{*}$ on $X$
such that
\roster
\myitem i.  $D^{*}-D$ is effective;
\myitem ii.  for all $e\in\Bbb Z_{>0}$ there exists a smooth variety $X'$
and a proper generically finite morphism $\pi\:X'\to X$ such that
$(\pi^{*}D)_\red$ has normal crossings and all components of $\pi^{*}D$
have multiplicity $\ge e$ (or zero);
\myitem iii.  $\pi\:X'\to X$ is unramified outside of $\pi^{-1}(\Supp D)$; and
\myitem iv.  $K(X')=K(X)(\root e\of{f_1},\dots,\root e\of{f_n})$ for some
$f_1,\dots,f_n\in K(X)^{*}$.
\endroster
\endit

\demo{Proof}  Let $n=\dim X$, and let $D_1,\dots,D_n$ be effective divisors
on $X$ such that $D_i\sim D$ for all $i$, $D^{*}:=D+D_1+\dots+D_n$ has
normal crossings, and $D\cap D_1\cap\dots\allowmathbreak\cap D_n=\emptyset$.
(Such divisors exist by Bertini's theorem.)

Now pick $e\in\Bbb Z_{>0}$.  For $i=1,\dots,n$ let $\pi_i\:X_i\to X$
be as in Lemma \03.3, applied to the divisors $D_i\sim D$.  Let $\pi\:X'\to X$
be a desingularization of the normalization of $X$ in the compositum
$K(X_1)\dotsm K(X_n)$ such that $(\pi^{*}D^{*})_\red$ has normal crossings
and $X'$ dominates $X_1,\dots,X_n$.  We may also assume that $\pi$
is \'etale outside of $\pi^{-1}(\Supp D)$.  Pick a component $E$ of
$\pi^{*}D^{*}$.  There exists some $i$ such that $\pi(E)\subseteq D\cup D_i$
but $\pi(E)\nsubseteq D\cap D_i$; then the image of $E$ in $X_i$ is contained
in a component of $\pi^{-1}(D+D_1)$ of multiplicity $\ge e$.\qed
\enddemo

After replacing $D$ with $D^{*}$, we have the following situation:
$D$ has normal crossings, $X'$ is smooth, $D':=(\pi^{*}D)_\red$
has normal crossings, all components of $\pi^{*}D$ have
multiplicity $\ge e$, and $\pi\:X'\to X$ is unramified outside of $D'$.

The methods of \cite{K, Thm.~17} and \cite{A-K, Lemma~5.8} allow one to
construct such a $\pi$ that is finite, but this is not necessary for the
present argument.

The key step in the proof of Theorem \03.1 is the following lemma,
which gives a bound on $d_k(P')-d_k(P)$.

\lemma{\03.5}  Pick models $\Cal X$ and $\Cal X'$ for $X$ and $X'$,
respectively, for which the divisors $D$ and $D'$ extend to effective
divisors on $\Cal X$ and $\Cal X'$, respectively, and $\pi$ extends to
a morphism $\pi\:\Cal X'\to\Cal X$.  Suppose $S$ contains all places lying
over $e$, all places of bad reduction for $\Cal X$ and $\Cal X'$, all places
where $D$ and $D'$ have vertical components, and all places where $\pi$ is
ramified outside the support of $D'$.  Let $r\in\Bbb Z_{>0}$.  For points
$P\in X(\bar k)\setminus\Supp D$ of degree $\le r$ over $k$ and
points $P'\in\pi^{-1}(P)$, we have
$$d_k(P') - d_k(P)
  \le N^{(1)}(D,P) - N(D',P') + \frac1e\,h_D(P) + O(1)\;,\tag\03.5.1$$
where the constant in $O(1)$ depends only on $X$, $D$, $e$, $X'$, $\pi$, $r$,
and $S$.
\endit

\demo{Proof}  By (\01.1) we have
$$d_k(P') - d_k(P) = d_{k,S}(P') - d_{k,S}(P) + O(1)\;.\tag\03.5.2$$
Let $E=k(P)$ and $E'=k(P')$.  For places $w'\in M_{E'}$ lying over places
$w\in M_E$, let $e_{w'/w}$ denote the index of ramification of $w'$
over $w$.  By (\01.2), we then have
$$d_{k,S}(P') - d_{k,S}(P)
  = \sum_{w'\in M_{E'}\setminus S} \frac{e_{w'/w}-1}{[E':k]}\mu(K(w'))\;.$$
If the closure in $\Cal X'$ of $P'$ does not meet $D'$ at $w'$,
then $E'$ is unramified over $E$ at $w'$; if it does meet, then $e_{w'/w}\le e$
because $E'$ is generated over $E$ by $e^{\text{th}}$ roots of elements of $E$
(by condition (iii) in Lemma \03.4).  Thus
$$d_{k,S}(P') - d_{k,S}(P) \le (e-1)N^{(1)}(D',P')\;.$$
Combining this with (\03.5.2) then gives
$$d_k(P') - d_k(P) \le (e-1)N^{(1)}(D',P') + O(1)\;.\tag\03.5.3$$

Now consider the right-hand side of (\03.5.1).  Since all components
of $\pi^{*}D$ have multiplicity $\ge e$,
$$N^{(1)}(D,P) - N^{(1)}(D',P') \ge (e-1) N^{(1)}(D',P')\;.$$
Also,
$$\split N(D',P') - N^{(1)}(D',P') &\le N(D',P') \\
  &\le \frac1e N(D,P) \\
  &\le \frac1e h_D(P) + O(1)\;.\endsplit$$
Combining these two inequalities then gives
$$N^{(1)}(D,P) - N(D',P') + \frac1e\,h_D(P) \ge (e-1)N^{(1)}(D',P') + O(1)\;.$$
Comparing this with (\03.5.3) then gives (\03.5.1).\qed
\enddemo

To prove the theorem, it remains only to assemble the ingredients.

\demo{Proof of Theorem \03.1}  By Kodaira's lemma we may assume (after
adjusting $\epsilon$) that $\Cal A$ is ample.  Then
$$h_D \le c\,h_{\Cal A} + O(1)\tag\03.6$$
for some constant $c$ depending only on $X$, $D$, and $A$.
Pick $e\ge c/\epsilon$, and let $X'$ be a generically finite cover
of $X$ as in Lemma \03.4.  Let $D$ and $D'$ be as discussed following
that lemma, and enlarge $S$ so that the conditions of Lemma \03.5 hold.
This will ultimately give an inequality (\02.3.1) relative to this larger
set $S$, which trivially implies the same inequality for the original $S$.
Points $P\in X(\bar k)\setminus\Supp D$ of bounded degree lift to
points $P'\in X'(\bar k)$, also of bounded degree.  We now show that
Conjecture \02.1 applied to $X'$ and $D'$ implies Conjecture \02.3
for $D$ and $X$.  By the former conjecture (using (\02.2)), we have
$$N(D',P') + d_k(P)
  \ge h_{\Cal K'(D')}(P') - \epsilon'h_{\Cal A'}(P') - O(1)\;,\tag\03.7$$
provided $P'\notin Z'$, where $Z'$ is a proper Zariski-closed subset
(depending also on $\epsilon'$ and $\Cal A'$); here $\Cal A'$ is a big line
sheaf on $X'$.  We want to show that (\03.7) implies Conjecture \02.3; \ie,
$$N^{(1)}(D,P) + d_k(P)
  \ge h_{\Cal K(D)}(P) - \epsilon\,h_{\Cal A}(P) - O(1)\;.\tag\03.8$$

By Lemma \03.2, $\Cal K'(D')\otimes\pi^{*}(\Cal K(D))\spcheck$ has a global
section vanishing nowhere on $X'\setminus\Supp D'$; hence, by functoriality
of heights and positivity of heights relative to effective divisors,
$$h_{\Cal K'(D')}(P') \ge h_{\Cal K(D)}(P) - O(1)\tag\03.9$$
for all $P\in X'(\bar k)\setminus D'$.

By Lemma \03.5,
$$N(D',P') + d_k(P')
  \le N^{(1)}(D,P) + d_k(P) + \frac1e\,h_D(P) + O(1)\;.\tag\03.10$$
Let $\Cal A'=\pi^{*}\Cal A$; this is big on $X'$; also choose $\epsilon'>0$
such that $\epsilon'<\epsilon-c/e$.  By (\03.6) we then have
$$\frac1e\,h_D(P) + \epsilon'h_{\Cal A'}(P')
  \le \epsilon\,h_{\Cal A}(P) + O(1)\;.\tag\03.11$$
Thus, (\03.9)--(\03.11), combined with (\03.7), imply that (\03.8) holds if
$P$ lies outside the proper Zariski-closed subset $Z:=\pi(Z')$.\qed
\enddemo

\enddocument